\documentclass[10pt,twoside,leqno]{amsart}

\usepackage{amssymb}
\begin{document}
\setcounter{page}{1}


\vspace*{1.5cm}

\title[]{\large Comments on relaxed $(\gamma , r)$-cocoercive mappings}
\author{ Shahram Saeidi}
\date{}
\maketitle

\vspace*{-0.5cm}

\begin{center}
{\footnotesize \textit{Department of Mathematics,  University of Kurdistan,
  Sanandaj 416, Kurdistan, Iran\\ Fax: +98 8713239948\\
E-mail: sh.saeidi@uok.ac.ir}}
\end{center}

\bigskip

\noindent {\bf Abstract.} We show that the variational inequality
$VI(C,A)$ has a unique solution for a relaxed $(\gamma ,
r)$-cocoercive, $\mu$-Lipschitzian mapping $A: C\to H$ with
$r>\gamma \mu^2$, where $C$ is a nonempty closed convex subset of a
Hilbert space $H$. From this result, it can be derived that, for
example, the recent algorithms given in the references of this
paper, despite their becoming more complicated, are not general as
they should be.\\


\noindent

\section*{}
Let $H$ be a Hilbert space, whose inner product is denoted by
$\langle . \;,.\rangle$. Let $C$ be a nonempty closed convex subset
of $H$ and let $A: C\to H$ be a nonlinear map. Let $P_C$ be the
projection of $H$ onto the convex subset $C$. The classical
variational inequality which is denoted by $VI(C,A)$ is used to find
$u\in C$ such that $$\langle Au, v-u\rangle\geq 0$$ for all $v\in
C$. For a given $x\in H$, $u\in C$ satisfies the inequality
$$\langle x-u, u-y\rangle\geq 0, \;\;\forall \; y\in C,$$
if and only if $u=P_Cx$. It is known that projection operator $P_C$
is nonexpansive.

It is easy to see that $$u\in VI(C,A)\Longleftrightarrow
u=P_C(u-\lambda Au),\;\;\;(1)$$ where $\lambda>0$ is a constant.
This alternative equivalent formulation has played a significant
role in the studies of the variational inequalities and related
optimization problems.
Recall that:\\

\noindent {\bf (i)} $A$ is called $r$-strongly monotone, if for each
$x, y\in C$ we have $$\langle Ax-Ay, x-y\rangle\geq r\|x-y\|^2$$ for
a constant $r>0$.\\
\noindent {\bf (ii)} $A$ is called relaxed $(\gamma, r)$-cocoercive
if there exists two constants $\gamma>0$ and $r>0$ such that
$$\langle Ax-Ay, x-y\rangle\geq -\gamma \|Ax-Ay\|^2+r\|x-y\|^2.$$
This class of maps is more general than the class of strongly
monotone maps.

 In [1], Verma considered a relaxed $(\gamma,
r)$-cocoercive $\mu$-Lipschitz mapping to present approximate
solvability of a system of variational inequality problems. In
[2], the same author proved the following theorem:\\

\noindent {\bf {Theorem}}{ (Verma [2]).} Let $C$ be a nonempty
closed convex subset of $H$ and let $A: C\to H$ be $r$-strongly
monotone and $\mu$-Lipschitz. Suppose that $x^*, y^*\in C$ be chosen
such that
$$\left\lbrace
  \begin{array}{c l}
    x^*=P_C(y^*-\rho Ay^*)\;\;\;\;\;\;\;\;  \text{for }\rho>0,\\
   y^*=P_C(x^*-\eta Ax^*)\;\;\;\;\;\;\;\;  \text{for }\eta>0.
  \end{array}
\right.\;\;\;\;\;\;\;(2)$$ For arbitrary chosen initial points $x_0,
y_0\in C$, define $\{x_n\}$ and $\{y_n\}$ as
$$\left\lbrace
  \begin{array}{c l}
    x_{n+1}=(1-a_n)x_n+a_n P_C(y_n-\rho Ay_n)\\
   y_n=(1-b_n)x_n+b_n P_C(x_n-\eta Ax_n)
  \end{array}
\right.\;\;\;\;\;\;\;(3)$$ for all $n\geq 0$, where, $0\leq a_n$,
$b_n\leq 1$ and $\sum_{n=0}^{\infty}a_nb_n=\infty$. Then sequences
$\{x_n\}$ and $\{y_n\}$, respectively, converge to $x^*$ and $y^*$
for
$$0<\rho<\frac{2r}{\mu^2}\;\hbox{and}\; 0<\eta<\frac{2r}{\mu^2}.$$

Very recently, Noor [3], Noor and Huang [4-7] and Qin et al. [8]
have considered some iterative methods for finding a common element
of the set of the fixed points of nonexpansive mapping and the set
of the solution of a variational inequality $VI(C, A)$, where $A$ is
a relaxed $(\gamma , r)$-cocoercive $\mu$-Lipschitzian mapping of
$C$ into $H$ such that $r>\gamma \mu^2$.

Moreover, Gao and Guo [9] and Qin et al. [10, 11]  proposed
iterative algorithms to find a common element of the set of
solutions of variational inequalities for a relaxed $(\gamma ,
r)$-cocoercive $\mu$-Lipschitzian mapping of $C$ into $H$ such that
$r>\gamma \mu^2$, the set of solutions of an equilibrium problem and
also the set of the fixed points of nonexpansive mapping.

By proving the following proposition we give some important comments
on the above mentioned works.\\

\noindent {\bf {Proposition 1.}} Let $C$ be a nonempty closed convex
subset of $H$ and let $A: C\to H$ be a relaxed $(\gamma ,
r)$-cocoercive, $0<\mu$-Lipschitzian mapping  such that $r>\gamma
\mu^2$. Then  $VI(C, A)$ is singleton.\\

\noindent { \textit{Proof.}} Let $s$ be a real number such that
$$0<s<\frac{2(r-\gamma \mu^2)}{\mu^2}.$$
Then, for every $x, y\in C$ we have
$$\|P_C(I-sA)x-P_C(I-sA)y\|^2$$$$\leq
\|(I-sA)x-(I-sA)y\|^2=\|(x-y)-s(Ax-Ay)\|^2$$$$ =\|x-y\|^2-2s\langle
x-y, Ax-Ay \rangle+s^2\|Ax-Ay\|^2$$$$\leq \|x-y\|^2-2s(-\gamma
\|Ax-Ay\|^2+r\|x-y\|^2)+s^2\|Ax-Ay\|^2$$$$\leq \|x-y\|^2
+2s\mu^2\gamma\|x-y\|^2-2sr\|x-y\|^2+\mu^2 s^2\|x-y\|^2$$$$
=(1+2s\mu^2\gamma-2sr+\mu^2
s^2)\|x-y\|^2$$$$=(1-s\mu^2[\frac{2(r-\gamma\mu^2)}{\mu^2}-s])\|x-y\|^2.$$
Now, since $1-s\mu^2[\frac{2(r-\gamma\mu^2)}{\mu^2}-s]<1$, the
mapping $P_C(I-sA): C\to C$ is a contraction and Banach' s
Contraction Mapping Principle guarantees that it has a unique fixed
point $u$; i.e., $P_C(I-sA)u=u$, which is the unique solution of
$VI(C, A)$ by (1). $\Box$\\

Since  $r$-strongly monotone mappings are relaxed $(\gamma, r)$-cocoercive, we get the following.\\

\noindent {\bf {Proposition 2.}} Let $C$ be a nonempty closed convex
subset of $H$ and let $A: C\to H$ be an $0<r$-strongly monotone and
$0<\mu$-Lipschitzian mapping. Then $VI(C, A)$ is singleton.\\

The followings are our comments:\\

\noindent {\bf {Comment 1.}} In Verma' s theorem, mentioned above,
$A: C\to H$ is an $0<r$-strongly monotone, $0<\mu$-Lipschitzian
mapping. Therefore, according to Proposition 2,  $VI(C, A)$ is
singleton, e.g.,   $VI(C, A)=\{u\}$. Takeing $x^*=y^*=u$, and
considering (1),  the equations in (2) hold and indeed  the
sequences $\{x_n\}$ and $\{y_n\}$ in (3) converge to the unique
solution $u$ of $VI(C, A)$. Therefore,  $x^*, y^*$ in Verma' s
theorem are indeed equal  the unique solution of $VI(C, A)$; thus
the system of variational inequality problems that is considered by
Verma in [2] comprises only one problem.\\

\noindent {\bf {Comment  2.}} The iterative methods considered by
Noor [3], Noor and Huang [4-7] and Qin et al. [8] for finding the
common element of the set of the solution of the variational
inequality $VI(C, A)$, where $A$ is a relaxed $(\gamma ,
r)$-cocoercive $\mu$-Lipschitzian mapping of $C$ into $H$ such that
$r>\gamma \mu^2$, and the set of the fixed points of nonexpansive
mapping, can not be really considered as algorithms for finding
fixed points of nonexpansive mappings; because, according to
Proposition 1, $VI(C, A)$ is singleton. Therefore, they must be
considered as algorithms converging only to the
unique solution of $VI(C,A)$.\\

\noindent {\bf {Comment 3.}} The  iterative algorithms proposed by
Gao and Guo [9] and Qin et al. [10, 11] for finding a common element
of the set of solutions of variational inequalities for a relaxed
$(\gamma , r)$-cocoercive $\mu$-Lipschitzian mapping of $C$ into $H$
such that $r>\gamma \mu^2$, the set of solutions of an equilibrium
problem and set of the fixed points of nonexpansive mapping,  are
not general as they should be and despite their becoming more
complicated, they converge only to the unique solution of $VI(C,A)$.

\end{document}